\documentclass{article}
\newcommand{\sect}[1]{\section{#1}}

\newtheorem{Th}{Theorem}
\newtheorem{Cor}[Th]{Corollary}

\def\proof#1. {\par
                      \ifdim\lastskip<15pt
                      \removelastskip\penalty-200
                      \vskip15pt plus3pt minus3pt
                      \fi
                       {\def\a{#1}
                       \ifx\a\empty
                       {\noindent\bf Proof.}
                       \else
                       {\noindent\bf Proof of #1.}
                       \fi}\enspace}
\def\restr#1{\,\vrule\,\lower1.75ex\hbox{$#1$}}

\def\endproof{\hfill\hspace{-6pt}\rule[-14pt]{6pt}{6pt}
\vskip22pt plus3pt minus 3pt}

\def\be{\begin{equation}}
\def\ee{\end{equation}}
\def\bea{\begin{eqnarray}}
\def\eea{\end{eqnarray}}
\def\bean{\begin{eqnarray*}}
\def\eean{\end{eqnarray*}}

\def\a{\alpha}
\def\b{\beta}
\def\d{\delta}
\def\D{\Delta}
\def\e{\varepsilon}
\def\f{\varphi}
\def\F{\Phi}
\def\g{\gamma}
\def\G{\Gamma}
\def\i{\infty}
\def\k{\kappa}
\def\l{\lambda}

\def\o{\omega}
\def\O{\Omega}
\def\r{\rho}

\def\t{\theta}

\def\z{\zeta}
\def\c{{\rm C}}

\def\setm{\setminus}

\def\ov{\overline}

\def\c{{\rm C}}

\font\tenopen = cmbx10
\font\sevenopen = cmbx7
\font\fiveopen = cmbx5
\newfam\openfam
\def\open{\fam\openfam\tenopen}
\textfont\openfam = \tenopen
\scriptfont\openfam = \sevenopen
\scriptscriptfont\openfam = \fiveopen

\def\R{{\open R}}

\def\C{{\open C}}

\usepackage{graphicx}
\date{\it Dedicated to the memory of two outstanding
mathematicians, Andrei Aleksandrovich Gonchar
and Herbert Stahl}
\title{On a conjecture of Widom\footnote{AMS Classification 42C05, 31A15,
Keywords: Widom's theory, Chebyshev polynomials, supremum norm, Jordan arcs}}

\author{Vilmos Totik\thanks{Supported by the European Research
Council Advanced Grant No. 267055}\  \ \ and Peter Yuditskii\thanks{Supported by the Austrian Science Fund, project no: P22025-N18}}
\begin{document}
\maketitle

\begin{abstract} In 1969 Harold Widom published his seminal paper \cite{W}
which gave a complete description of orthogonal and
Chebyshev polynomials on a system of smooth Jordan curves.
When there were Jordan arcs present the theory of orthogonal
polynomials turned out to be just the same,
but for Chebyshev polynomials Widom's
approach proved only an upper estimate,
which he conjectured to be the correct
asymptotic behavior. In this note we make some clarifications
which will show that
the situation is more complicated.
\end{abstract}

\sect{Widom's problem for Chebyshev polynomials} This paper uses some basic facts from
logarithmic potential theory, see the books
\cite{Gardiner}, \cite{Garnett} or \cite{Ransford}
 for the concepts used.

Let $E$ be the union of finitely many disjoint
Jordan curves or arcs $E_k$  that lie in the exterior of each
other and that satisfy some smoothness condition.
Recall that a Jordan curve is a homeomorphic image of
the unit circle, while a Jordan arc is a homeomorphic image
of a segment.
The Chebyshev polynomial of degree $n$ associated
with $E$ is the unique polynomial $T_n(z)=z^n+\cdots$
which minimizes the supremum norm
\[\|T_n\|_E=\sup_{z\in E}|T_n(z)|.\]
If the minimal norm is denoted by $M_n$, then
it is known (see e.g. \cite[Theorem 5.5.4]{Ransford}) that
\be M_n\ge \c(E)^n,\qquad n=1,2,\ldots,\label{also1}\ee
where $\c(E)$ denotes logarithmic capacity, and $M_n^{1/n}\to \c(E)$
as $n\to\i$ (\cite[Corollary 5.5.5]{Ransford}).
It is a delicate problem how close $M_n$ can get to
the theoretical lower bound $C(E)^n$; the questions
we are dealing with in this paper are also connected
with that problem.

It is easy to see that if $E$ is the unit circle then
$M_n=1=\c(E)^n$, and in general, Faber proved \cite{Faber} that
for a single Jordan curve $M_n\sim \c(E)^n$,
where $\sim$ means that the ratio of the two sides
tends to 1. The original problem that was considered by Chebyshev was
for $E=[-1,1]$, in which case $M_n=2\cdot 2^{-n}=2\cdot \c(E)^n$,
twice $\c(E)^n$.

If $\r$ is a nonnegative weight function on $E$ then
one can similarly define the weighted Chebyshev polynomials $T_{n,\r}$
and the weighted Chebyshev numbers $M_{n,\r}$. If we use the square of
the $L^2(\r)$-norm instead of $L^\i(\r)$
(i.e. consider minimizing
\[\int_E |Q_n(z)|^2\r(z)d|z|\]
instead of
\[\sup_E |Q_n(z)|\r(z)),\]
 then one obtains the
quantities $m_{n,\r}$ and the extremal polynomials
$P_{n,\r}(z)=z^n+\cdots$, and here $P_{n,\r}$
turns out to be the $n$-th monic orthogonal polynomial with respect
to $\r$ and $P_n(z)/\sqrt{m_{n,\r}}$ is the orthonormal one.
With these notations (under suitable conditions on
$E$ and $\r$) in the paper \cite{W} Harold
Widom gave a full description of the
quantities $M_{n,\r}$, $m_{n,\r}$, and also
determined the precise behavior of the Chebyshev and
orthogonal polynomials themselves away from the set $E$,
provided $E$ consists of Jordan curves. The description was
in turn of some associated Green and Neumann functions.
When some of
the components of $E$ were Jordan arcs then
he also proved the same behavior for the
orthogonal polynomials $P_{n,\r}$ and for $m_{n,\r}$,
but his theory was not complete in that case for
the Chebyshev polynomials $T_{n,\r}$ and the
Chebyshev numbers $M_{n,\r}$. He wrote in discussing the interval case:
``Thus $M_n$ is asymptotically twice
as large for an interval as for a closed curve of
the same capacity. We conjecture that this is
true generally: that is, if at least one of the $E_k$
is an arc then the asymptotic formula
for $M_{n,\r}$ given in Theorem 8.3 must be
multiplied by 2.\ $\ldots$\  Unfortunately we cannot prove these
statements and so they are nothing but conjectures".
Widom himself showed that his conjecture
is true if $E$ lies on the real line,
i.e. it consists of a finite number
of intervals,
and he also showed that the asymptotics for $M_{n,\r}$
is {\it at most twice} as large as the
asymptotics in the curve case.

In particular, for $\r\equiv 1$ the conjecture would
imply that if $E$ has an arc-component, then
\be \liminf_{n\to\i}\frac{M_n}{\c(E)^n}\ge 2,\label{2}\ee
since (\ref{also1})
is true for any set $E$.

Widom's paper had a huge impact on the theory
of extremal polynomials, in particularly on
the theory of orthogonal polynomials. It is
impossible to list all further contributions,
for orientation see
the papers \cite{Aptekarev1}--\cite{Aptekarev2},
\cite{Kal},
\cite{Peh1}--\cite{Peh5}. See also
\cite{Schief} for a lower bound for
the Chebyshev constants for sets on the real line,
and the paper \cite{Sodin} for the various connections
of the Chebyshev problem.

The aim of this note is to make some simple clarifications
in connection with Widom's conjecture. Strictly speaking the conjecture
is not true in the stated form, but we shall
see that Widom was absolutely right that arcs
change the asymptotics when compared to asymptotics on curves.
Originally the authors had some ideas indicating
that the situation was more complex
than how Widom conjectured, but then they realized
that more than what they wanted to say can be
deduced rather easily from Widom's work itself,
so this note follows the setup and the
arguments in \cite{W} very closely.

For the case when $E$ consisted purely of Jordan curves
the asymptotics of Widom was in the form
\be M_{n,\r}\sim C(E)^n\mu(\r,\G_n)\label{w1}\ee
with some rather explicitly given quantity
$\mu(\r,\G_n)$; see the next section. In the general case
when $E$ may have curve and arc components the following holds.
\begin{Th}\label{Thmain} If there is at least one Jordan
curve in $E$, then
\be \limsup_{n\to\i}\frac{M_{n,\r}}{\c(E)^n\mu(\r,\G_n)}\le \t<2,\label{hhh}\ee
where $\t$ depends only on $E$.
\end{Th}

This of course disproves (\ref{2}), however (\ref{2}) is partially true:
it was proved in \cite[Theorem 1]{TotikCG} that if
a general compact $E$ contains
an arc on its outer boundary,
then there is a $\b>0$ such that
\[\liminf_{n\to\i}\frac{M_n}{\c(E)^n}\ge 1+\b.\]
In general one cannot say much more than that.
In fact, it was proven by Thiran and Detaille \cite[Section 5]{TD} that if
$E$ is a subarc on the unit circle of central angle
$2\a$, then
\[T_n\sim C(E)^n 2\cos ^2\a/4.\]
The factor $2\cos ^2\a/4$ on the right
is always smaller than 2 and is as close
to 1 as one wishes if $\a$ is close to $\pi$.

As for asymptotics for $M_{n,\r}$, we shall prove in
Section \ref{Sec3} the following. Let $\r^*$ be equal
to $\r$ on the curve components of $E$ and equal to $2\r$ on
the arc components. With this (\ref{hhh}) is a consequence
of
\be \limsup_{n\to\i}\frac{M_{n,\r}}{\c(E)^n\mu(\r^*,\G_n)}\le 1,\label{hhh1}\ee
to be proven in the next section (see (\ref{3})). The next theorem
shows that this estimate is exact when the set is symmetric
with respect to the real line.

\begin{Th}\label{Thmain1} If $E$ consists of real intervals
 and of Jordan curves symmetric with respect
to the real line, then
\[M_{n,\r}\sim\c(E)^n\mu(\r^*,\G_n).\]
\end{Th}

It is known that if $g(z,\i)$ is Green's function of
the outer component of $E$ with pole at infinity, then
$g$ has $p-1$ critical points which we denote as
$z_1^*,\ldots,z_{p-1}^*$. Furthermore, let
$\nu_E$ denote the equilibrium measure of $E$, and let
$E_{\rm arc}$ be
the union of the arc-components of $E$.

Theorem \ref{Thmain1} gives the following bound for the Chebyshev
constants.
\begin{Cor}\label{cor} Under the conditions of Theorem \ref{Thmain1} the limit points
of the sequence $\{M_n/C(E)^n\}$ lie in the interval
\be \left[ 2^{\nu_E(E_{\rm arc})},2^{\nu_E(E_{\rm arc})}\exp\left\{\sum_{j=1}^{p-1}g(z_j^*)
\right\}\right].\label{b}\ee
Furthemore, if
\[\nu_E(E_1),\ldots,\nu_E(E_p)\]
are rationally independent, then the set of limit points
is precisely the interval (\ref{b}).\end{Cor}

In particular, if $E=E_{\rm arc}$, i.e. $E$ is a subset
of the real line, then
\[\liminf_{n\to\i}\frac{M_n}{\c(E)^n}\ge 2,\]
as was proved by Widom, see also
\cite{Schief}.

We note that in the case when $E$ is consisting of
Jordan curves we have $\nu_E(E_{\rm arc})=0$, and so the corresponding interval
is
\[ \left[ 1,\exp\left\{\sum_{j=1}^{p-1}g(z_j^*)
\right\}\right],\]
see \cite[Theorem 8.4]{W}.

In the next two sections we show how to prove Theorems \ref{Thmain}
and \ref{Thmain1} using the setup and reasonings of Widom's paper \cite{W}.
In the last section we give an explicit formula (\ref{ex}) for the asymptotics $M_{n,1}$ in the elliptic case, that is, in the case that the boundary of $\Omega$ consists of  a real interval $[\alpha,\beta]$ and a symmetric curve. The ratio $M_{n,1}/C(E)^n$ behaves in $n$ 
as an almost periodic (or periodic) function, 
depending on  the modulus of the domain and 
on the harmonic measure of the interval 
evaluated at infinity (which is the same as
the mass of the equilibrium measure carried by the
interval). 
Let us mention that generally research in this direction was started in \cite{A1}.

\sect{Widom's theory and Theorem \ref{Thmain}}
We shall need to briefly describe Widom's paper
\cite{W}.

Let $E=\cup_{k=1}^p E_k$ be a finite family of Jordan curves and arcs
lying exterior to one another.
The smoothness assumptions on $E$ we take the assumptions of Widom's paper
\cite{W}, $C^{2+}$ will certainly suffice. We shall denote by
$E_{\rm arc}$ the union of the arc components of $E$.
Let further $\r$ be a weight on E, of which we assume for simplicity
that it is positive and satisfies a Lipshitz condition, i.e. it is
of class $C^+$.

The weighted Chebyshev numbers with respect to $\r$  will be denoted by
$M_{n,\r}$, i.e.
\[ M_{n,\r}=\inf \|\r(\z)(\z^n+\cdots)\|_E,\]
where $\|\cdot\|_E$ denotes the supremum norm on $E$ and
where the infimum is taken for all monic polynomial
of degree $n$.

$\O$ denotes the outer domain, i.e. the unbounded component of $\ov \C\setm E$,
$g(z,w)$ is its Green's function with pole at $w$, and
$\F(z,w)=\exp(g(z,w)+i\tilde g(z,w))$ is the so called complex
Green's function with the normalization that $\F(\i,w)>0$ if $w$ is finite
and $\F(z,\i)/z\to c>0$ as $z\to \i$ if $w=\i$. Here $c=1/C(E)$,
the reciprocal of the logarithmic capacity $C(E)$ of $E$.
The $\F$ is a multi-valued analytic function,
but $|\F|$ is single-valued. Consider
 multi-valued bounded functions $F$ in $\O$ for which $|F|$ is single-valued
 and which have only finitely many zeros in $\O$. Each such $F$ has boundary values on $E$ almost everywhere.
If we
 define
 \[\g_k=\g_k(F)=\frac1{2\pi}\D_{E_k}\arg F,\qquad k=1,\ldots,p,\]
 where the total change $\D_{E_k}\arg F$ of the argument of $F$
 around $E_k$ is taken on some positively oriented
 curve in $\O$ lying close
 to $E_k$, then $\G(F)=(\g_1,\ldots,\g_p)$ is called the class
 of $F$.
The class of $\F^{-n}$ is denoted by $\G_n$, $n=1,2,\ldots$.

For a given class $\G=(\g_1,\ldots,\g_p)$ let
$\mu(\r,\G)$ be the minimum of the norms $\sup_E \r |F|$,
where $F$ runs through all functions in the class $\G$
with the property that $F(\i)=1$. There is a unique
extremal function minimizing this norm, and it is of the form (see
\cite[Theorem 5.4]{W})
\be F_{\r,\G}(z)=\mu(\r,\G)R^{-1}(z)\prod \F(z,z_j)^{-1},\label{extremal}\ee
where ${z_1,\ldots z_q}$ are some points in $\O$, their number
$q$ is at most $p-1$, and
where $R(z)=R_\r(z)$ is the outer function in $\O$ with boundary values $\r$,
i.e. $R$ is the multi-valued analytic function in $\O$ which is positive
at $\i$, and for which  $\log|R(z)|$ is single-valued and has boundary values $\r$ on $E$.
The choice of the points $z_1,\ldots,z_q$ is such that
$F$ is of class $\G$, and with them we have for the extremal constant
the expression
\be \mu(\r,\G)=|R(\i)|\exp\left\{\sum_{j=1}^q g(z_j)\right\}.\label{muform}\ee

We shall also need the harmonic measures $\o_k(z)$ in
$\O$ associated with $E_k$, i.e. $\o_k(t)$ is
the harmonic function in $\O$ which has boundary value 1 on $E_k$
and 0 on the rest of $E$. If $\tilde \o_k$ denotes its
analytic conjugate with some normalization, then $\exp(\o_k(z)+i\tilde \o_k(z))$
is again a multi-valued analytic function in $\O$
for which its absolute value is single-valued.
When doing asymptotics, it is desirable not to have the $x_j$'s from
(\ref{extremal}) lie too close to $E$, and in that case
Widom changed $F_{\r,\G}$ from (\ref{extremal}) to (see
\cite[p. 215-216]{W})
\be F_\G(z)=
F_{(\r,\G)}=\mu(\r,\G)V_{\G'}^{-1}(\i)R^{-1}(z)\prod_1 \F(z,z_j)^{-1}V_{\G'}(z),\label{extremal1}\ee
where in $\prod_1$ those $\F(z,z_j)^{-1}$ are kept for which
$z_j$ are
of distance $\ge \d$ from $E$ for some small fixed $\d$, and
\[V_{\G'}(z)=\exp\left\{\sum_{k=1}^p\l_k(\o_k(z)+i\tilde \o_k(z))\right\}\]
with some appropriate $\l_k$ that ensures that $F_\G$ is still of class
$\G$. This can be done for all $\G$ uniformly, and although this
$F_\G$ is no longer extremal for $\mu(\r,\G)$,
it is a nice smooth function on $E$ (uniformly in $\G$),
 and the norm $\sup_E\r |F_{\G}|$ is close
to $\mu(\r,\G)$ (again uniformly in $\G$). In fact, $\l_k$
are small if $\d>0$ is small, and
$\r(\z)|F_\G(\z)|/\mu(\r,\G)$ is uniformly as close to 1
as we wish if $\d>0$ is sufficiently small.

\bigskip

Define the weight $\r^*$ equal to $2\r$ on $E_{\rm arc}$ and equal to $\r$ on the rest
of $E$ (i.e. on the curve components of $E$), and consider the
functions just introduced, but for $\r^*$ rather than for $\r$. In particular,
consider $F_{(\r^*,\G)}$ from (\ref{extremal1})
with $\r$ replaced by $\r^*$. In what follows let
$\G_n$ be the class of $\F^{-n}(z)$, and consider, as Widom did,
\[Q(z)=\frac{1}{2\pi i}\int_C R_{\r^*}^{-1}(z)\prod_1 \F(z,z_j)^{-1}V_{\G_n'}(z) \F(\z)^n
\frac{d\z}{\z-z},\]
where $C$ is a large circle
about the origin (described
once counterclockwise) containing $E$ and $z$ in its interior.
We emphasize that the integrand is single-valued in
$\O$, and in the expression in the integrand in front of
$\F(\z)^n$ we have the function
from (\ref{extremal1}) (modulo some constants) made for
$\r^*$ and for the class $\G_n$. Now this $Q$ is a polynomial of degree
$n$, for which Widom proved in  Theorem 11.4 (when $\r$ is replaced by $\r^*$)
 the asymptotic formula as $n\to\i$
(see middle of \cite[p. 210]{W})
\[Q(\z)=B(\z)+o(1),\qquad \z\in E,\]
where $B(\z)$ is the limiting value of
\be R^{-1}_{\r^*}(z)\prod_1 \F(z,z_j)^{-1}V_{\G'}(z)\F(z)^n\label{pr1}\ee
on the closed curves of $E$, and the sum of the two limiting values
on the arcs of $E$. Here the second  and last factors  are of absolute
value 1 on $E$ while the first one is $1/\r^*(\z)$ and the third one is
as close to 1 as we wish, say $|V_{\G'}(\z)|\le e^{2\e}$ for any given
$\e>0$ provided the $\d>0$ above is sufficiently small
(and fixed for all $n$). This gives
(see \cite[p. 210]{W})
\[\sup_{E\setm E_{\rm arc}} |Q(\z)|\r(\z)\le e^{2\e}\]
and
\[\sup_{E_{\rm arc}} |Q(\z)|\r(\z)\le e^{2\e}2\sup_{E_{\rm arc}}\frac{1}{|R_{\r^*}(\z)|}\r(\z)
=e^{2\e}2\sup_{E_{\rm arc}}\frac{1}{2\r(\z)}\r(\z)\le e^{2\e}.\]
The absolute value of the leading coefficient is at least (see \cite[p. 210 and the proof
of Theorem 8.3]{W})
\[e^{-\e}C(E)^{-n}\mu(\r^*,\G_n),\]
and hence
\be M_{n,\r}\le e^{3\e} C(E)^n\mu(\r^*,\G_n).\label{3}\ee

Let $\t=\mu(\tau,0)$, where $\tau(\z)=2$ on $E_{\rm arc}$ and $\tau(\z)=1$ on $E\setm E_{\rm arc}$, and
where in $\mu(\tau,0)$ we take the 0-class of analytic functions in $\O$.
In other words,
$\t$ is the infimum of the norms $\sup_E \nu(\z) |h(\z)|$ for all $h$
which is bounded and analytic in $\O$ and equals 1 at infinity.
Then $\t<2$ provided there is
at least one curve component of $E$. This follows from the
fact that clearly $\mu(\tau,0)\le 2$ if we use
$h(z)\equiv 1$ as a test function, and that function
is not extremal, since for the extremal function (\ref{extremal})
the product
$\nu(\z)|h(\z)|$ is constant on $E$ (see \cite[Theorem 5.4]{W}). Now  the
definition of $\mu(\r,\G_n)$ implies that
$\mu(\r^*,\G_n)\le \t\mu(\r,\G_n)$, and hence (\ref{hhh})
 follows from (\ref{3}) (because $\e>0$ is arbitrary).\endproof

\sect{Widom's theory and Theorem \ref{Thmain1}}\label{Sec3}
In this section we assume, as in Theorem \ref{Thmain1},
that $E$ consists of intervals on the real line and of
Jordan curves that are symmetric with respect to the real line.
We show that a simple modification of some of Widom's argument
gives the asymptotic formula in Theorem \ref{Thmain1}.

The upper estimate is in (\ref{hhh1}) (see (\ref{3})), so we shall only deal with
the lower estimate of $M_{n,\r}$.

Let $z_1=z_{1,n},\ldots,z_{q_n,n}$ be the points from (\ref{extremal})
for $\G=\G_n$,
both their number $q=q_n$ and the points themselves depend on
$n$, but we shall suppress this dependence.

The intersection $E\cap \R$ consist of some intervals
$[\a_k,\b_k]$, $\a_1<\b_1<\a_2<\cdots<\a_p<\b_p$.
We call $(\b_k,\a_{k+1})$ the contiguous intervals to $E\cap \R$.
First we show that the points $z_j$ belong to the contiguous
intervals, and each contiguous interval contains at most
one of them. To this end note that $\O$ can be mapped
by some conformal map $\f$
onto some set $\ov\C\setm \cup_{k=1}^p [\a_k',\b_k']$
such that $(\b_k,\a_{k+1})$ is mapped into $(\b_k',\a_{k+1}')$.
Indeed, just take a conformal map of $\O\cap \C_+$ (where $\C_+$ is the upper
half plane) onto $\C_+$, and extend it across $\R\setm \cup_{k=1}^p [\a_k,\b_k]$
by Schwarz reflection. Now the extremal problem described
in the preceding section is conformally invariant,
so $\f(z_i)$ are mapped into the corresponding $z_j'$'s for the set
$E'=\cup_{k=1}^p [\a_k',\b_k']$ and for the function
$\f^{-1}(\r)$ keeping the same class $\G$. But simple variation argument
(see \cite[p. 211]{W}) shows that if $E'$ lies on the real line
then the correspond $z_j'$ lie in the contiguous intervals
to $E'$ and each of these intervals may contain at most one
of the $z_j'$, which proves the claim above.

Let
\be H(z)=R^{-1}_{\r^*}(z)\prod_1 \F(z,z_j)^{-1}V_{\G'}(z)\label{pr11}\ee
be the expression in (\ref{pr1}) in front of $\F(z)^n$.
Note that in (\ref{pr11}) not all $z_j$ appear, only those that are
of distance $\ge \d$ from $E$.
Let $l_k=1$ if the interval $(\b_k,\a_{k+1})$ contains
a $z_j$ appearing in (\ref{pr11}), and let otherwise $l_k=0$,
$1\le k\le p-1$.

\def\arg{{\rm arg}\ }
Now consider the proof of \cite[Theorem 11.5, pp. 212-214]{W}.
As there, we cut $\O$ along all contiguous intervals
and also along $(\b_p,\i)$, and consider
the argument of $H_\pm(\z)\F_\pm (\z)^n$.
Take the conjugate functions so that the argument of $H(z)\F^n(z)$
is
0 at $\a_1$. Then  it is 0 on all $(-\i,\a_1)$, and  $H(z)\F^n(z)$
is symmetric
with respect to the real line.
If $E_1$ is an interval ($[\a_1,\b_1]$), then
we get exactly as on  p. 213 of \cite{W} from the single-valuedness of
$H(z)\F(z)^n$ that
\be  \arg H_-(\b_1)\F_- ^n(\b_1)=m_1\pi\label{b1}\ee
\be  \arg H_+(\b_1)\F_+ ^n(\b_1)=-m_1\pi\label{b2}\ee
with some integer $m_1$. On the other hand,
if $E_1$ is a Jordan curve, then again the symmetry of
$E$ and the single-valuedness of $H(z)\F(z)^n$ imply
(\ref{b1})--(\ref{b2}). The number $m_1$ is positive
for large $n$ (actually tends to infinity as $n$ tends to infinity).
This follows from the Cauchy-Riemann equations on $E_1$ (more precisely
on some smooth curve lying close to $E_1$) and from the fact
that the normal derivative of $\log|H(z)|$ in the direction
of the inner normal to $\O$ is bounded, while the normal
derivative of $\log|\F(z)|$ is positive on $E_1$.

If $l_1=0$, then
the argument at $\a_2$ is the same, but if $l_1=1$,
then
\[ \arg H_-(\a_2)\F_- ^n(\a_2)=(m_1+1)\pi\]
\[ \arg H_+(\a_2)\F_+ ^n(\a_2)=-(m_1+1)\pi.\]

Proceeding this way, we get exacly as Widom that
\bean  \arg H_\pm(\b_k)\F_\pm ^n(\b_k)&=&\mp \left(\sum_{r=1}^k m_r +\sum_{r=1}^{k-1} l_k\pi\label{b11}\right)\pi\\
\arg H_\pm (\a_k)\F_\pm ^n(\a_k)&=&\mp \left(\sum_{r=1}^{k-1} m_r +\sum_{r=1}^{k-1} l_k\pi\right)\pi\nonumber\eean
with some positive integers $m_1,\ldots,m_p$.

On the other hand, it is also true (see p. 214 in \cite{W}) that on the arc components
of $E$
\be |H_\pm(\z)\F_\pm(\z)^n|\r^*(\z)=|H_\pm(\z)|2\r(\z)=1+O(\e),\label{kk}\ee
while on the curve components we have similarly
\be |H(\z)\F(\z)^n|\r^*(\z)=|H(\z)|\r(\z)=1+O(\e)\label{kk*}\ee
(where $\e>0$ depends on the $\delta$ that we used in
selecting or deleting the terms $\F(z,z_j)$ in (\ref{pr11})).
At the same time the polynomial
\[Q(z)=\int_C\frac{H(\xi)\F(\xi)^n}{\xi-z}d\xi\]
 constructed from $H(z)\F(z)^n$ satisfies (see \cite[Lemma 11.2]{W})
\[Q(\z)=H_+(\z)\F_+ (\z)^n+H_-(\z)\F_- (\z)^n+o(1)\]
on every arc-component of $E$ and
\be Q(\z)=H(\z)\F(\z)^n+o(1)\label{yu}\ee
on every curve-component of $E$.
Hence, if
\[\psi(\z)=\arg H_-(\z)\F_- (\z)^n,\]
then on any arc-component (cf. also (\ref{kk}))
\be Q(\z)\r(z)=\frac{1}{2}Q(\z)\r^*(z)=\frac{1}{2}(2\cos\psi(\z)+O(\e))=
\cos\psi(\z)+O(\e),\label{uu}\ee
while on the the curve-components (\ref{yu}) is true,
which gives, in view of (\ref{kk*}) that
\be |Q(\z)|\r(\z)=1+O(\e)\label{7}\ee
there. We also get from (\ref{yu}) and from the argument
principle that if $E_k$ is a curve-component
of $E$, then $Q$ has $m_k$ zeros inside $E_k$.
 Since on any arc-component $[\a_k,\b_k]$ the function $\psi$ is changing
from some $a\pi$ with some integer $a$
to $a\pi+m_k\pi$, we get from (\ref{uu}) that
there are $m_1+1$ points $\a_k\le x_{1,k}<\cdots<x_{m_k+1,k}\le \b_k$ where
$Q(\z)\r(\z)$ alternatively takes the values $\pm (1-\eta)$,
where $\eta$ is some small number that can
be as small as we wish if $\e>0$ is sufficiently small.
In addition, if $l_k=1$, then there is
an additional sign change along $(\b_k,\a_{k+1})$, i.e.
$Q(\b_k)\r(\b_k)$ and $Q(\a_{k+1})\r(\a_{k+1})$ are of opposite sign
and are close to 1, say $\ge 1-\eta$ in absolute value.
We may also assume that the $O(\e)$ in (\ref{7}) is
$\le \eta$ in absolute value.

Let $L=\sum_{k=1}^{p-1}l_k$ be the total number of the $z_j$'s
in (\ref{pr11}), which is the number of zeros of
$H(z)$ in $\O$. Since the change of the argument
of $H(z)\F^n(z)$ around a large circle is $2n\pi$, it
follows from the argument principle that the
total change of the argument of $H(z)\F^n(z)$
around $E$ is $2(n-L)\pi$. On the other hand, we have
calculated the total change of the argument
around $E$ to be $2\pi(m_1+\cdots +m_p)$,
hence $n=L+\sum_k m_k$.

Now these easily imply that if $P$ is any $n$-th degree
polynomial with the same leading coefficient as $Q$, then
\[\max_E|P(\z)|\r(\z)\ge 1-\eta.\]
Indeed, otherwise $Q-P$ would have alternating signs at the
$x_{i,k}$'s, giving $m_k$ zeros on every arc-component $[\a_k,\b_k]$.
With the same reasoning $Q-P$ has a zero on every
contiguous interval $(\b_k,\a_{k+1})$ for which $l_k\not=0$.
By (\ref{7}) and by the indirect assumption $\r(\z)|Q(\z)|<1-\eta$
we get from  Rouche's theorem that
$Q-P$ has the same number of zeros inside any curve-component
$E_k$ as $Q$ has there, i.e. $m_k$.
Thus, altogether we would get $\sum_k m_k +L=n$ zeros
for $Q-P$, which is impossible
since $Q-P$ is of degree $<n$.

As a consequence, we obtain that if $Q(z)=\k_nz^n+\cdots$,
then
\be M_{n,\r}\ge \frac{1-\eta}{|\k_n|}.\label{mm}\ee
By the definition of the extremal quantity $\mu(\r^*,\G_n)$ we have
\[ \sup_E \r^*(\z)|H(\z)/H(\i)|\ge \mu(\r^*,\G_n),\]
and since $\r^*(\z)|H(\z)|=1+O(\e)$ (see (\ref{kk})--(\ref{kk*})), it follows that
\[ \frac{1}{|H(\i)|}\ge e^{-\e}\mu(\r^*,\G_n)\]
provided $\d>0$ is sufficiently small. Therefore,\[\frac{1}{|\k_n|}=\frac{1}{|H(\i)|C(E)^{-n}}\ge C(E)^ne^{-\e}\mu(\r^*,\G_n),\]
and the lower estimate
\[\liminf\frac{M_{n,\r}}{C(E)^n\mu(\r^*,\G_n)}\ge 1\]
follows from (\ref{mm}).\endproof
\bigskip

Corollary \ref{cor} follows easily. Indeed,
(\ref{muform}) shows that
\be R_{\r^*}(\i)\le \mu(\r^*,\G_n)\le R_{\r^*}(\i)\exp\left\{\sum_{j=1}^{p-1}g(z_j^*)
\right\}.\label{vb}\ee
When $\r=1$ we have
\[|R_{\r^*}(z)|=\exp\{  \o_{E_{\rm arc}}(z)\ln 2\}\]
where $\o_{E_{\rm arc}}$ is the harmonic measure
in $\O$ corresponding to $E_{\rm arc}$, and it
is well known (see e.g. \cite[Theorem 4.3.14]{Ransford})
that
\[\o_{E_{\rm arc}}(\i)=\nu_E(E_{\rm arc})\]
where $\nu_E$ is the equilibrium measure of $E$.
Now Corollary \ref{cor} is a consequence of (\ref{vb}) and Theorem \ref{Thmain1}.

For the last statement in the corollary follow the
proof of \cite[Theorem 8.4]{W}.\endproof

\sect{Elliptic case}

In this section we give an explicit formula for the asymptotics of $M_{n,1}$ assuming that $E$ consists of an interval $[\alpha_1,\beta_1]$ and a symmetric Jordan curve $E_2$. In this case the domain $\Omega$ is conformally equivalent to an annulus $\{z:\ r_1<|z|<r_2\}$.  The ratio $r_2/r_1$ is a {\em conformal invariant of the domain} and the expression $\frac 1{2\pi}\ln r_2/r_1$ is called the {\em modulus} of $\Omega$ (it is the extremal length of curves in $\Omega$ that connect the two boundary components of $\Omega$). Following \cite[\textsection 55]{AKHef} we use the notation
$$
\tau:=\frac{i}{\pi }\ln\frac{r_2}{r_1}=2i\, {\rm mod}(\Omega).
$$

\begin{Th}
Let $\omega(\infty)$ be the harmonic measure of the interval $[\alpha_1,\beta_1]$ in $\Omega$ evaluated at infinity.
Then
\begin{equation}\label{ex}
M_{n,1}\sim C(E)^n 2^{\omega(\infty)}
\left|
\frac
{\vartheta_0\left(\frac{\{n\o(\infty)+|\tau'|\frac{\ln 2}{\pi} \}+\o(\infty)}{2}\left|\right. \tau' \right)}
{\vartheta_0\left(\frac{\{n\o(\infty)+|\tau'|\frac{\ln 2}{\pi} \}-\o(\infty)}{2}\left|\right. \tau' \right)}
\right|,
\end{equation}
where $\tau'=-1/\tau$, $\{x\}$  denotes the fractional part of a real number $x$, and
$$
\vartheta_0(t |\tau')=1-2h\cos 2\pi t+2h^4\cos 4\pi t-2h^9\cos 6\pi t+\dots,
\quad h=e^{\pi i\tau' },
$$
is the theta-function.
\end{Th}
Recall that here $\o(\i)=\nu_E([\a_1,\b_1])$, where
$\nu_E$ is the equilibrium measure.
\noindent
\proof. Using the conformal map $\f$ from Section \ref{Sec3} that maps $\O$
onto some $\overline{\C}\setm ([\a_1',\b_1']\cup [\a_2',\b_2'])$,
 we have the conformal mapping
$$
u=u(z)=\int_{\beta'_2}^{\f(z)}\frac{d\xi}{\sqrt{(\xi-\alpha'_1)(\xi-\beta'_1)(\xi-\alpha'_2)(\xi-\beta'_2)}}
$$
of $\Omega$ (which we cut along the contiguous interval $[\beta_1,\alpha_2]$) onto the rectangle with the vertices $-iK',K-iK', K+iK', iK'$,
$$
K=u(\alpha_1), \quad K+iK'=u(\beta_1).
$$
Therefore
$$
z\mapsto e^{-\frac{u(z)}{K'}\pi}
$$
maps conformally $\Omega$ onto the annulus $e^{-K\pi /K'}<|w|<1$, in particular $\tau=i\frac{K}{K'}$.

In this notations for the harmonic measure $\omega(z)=\omega(z,[\alpha_1,\beta_1])$ we have
$$
\omega(z)=\frac{1}{K}\Re u(z),
$$
and therefore
$$
\frac 1{2\pi}\Delta_{E_1}\arg \frac{1}{\Phi}=\omega(\infty)=\frac 1 K\Re u(\infty), \quad
\frac 1 {2\pi}\Delta_{E_1}\arg \frac{1}{\Phi(z,z_1)}=\omega(z_1)=\frac 1 K\Re u(z_1).
$$

Further,
 $R(z)=\exp(\ln 2 \frac{u(z)}{K})$. That is,
 \begin{equation}\label{ex7}
\Delta_{E_1}\arg R=\frac{2 K'}{K}\ln 2 \quad {\rm and}\quad R(\infty)=2^{\omega(\infty)}.
\end{equation}
Since  the product $R^{-1}(z)\Phi(z,z_1)^{-1}\Phi(z)^n$ is single-valued, we obtain the following expression for the real part of $u(z_1)$
$$
\frac 1 K\Re u(z_1)=\omega(z_1)=\frac{K'}{\pi K}\ln 2+n\omega(\infty)\ {\rm mod}\, 1.
$$
Since $0\le \omega(z_1)\le 1$ we get
$$
\frac 1 K\Re u(z_1)=\left\{\frac{K'}{\pi K}\ln 2+n\omega(\infty)\right\},
$$
and therefore
\begin{equation}\label{ex2}
u(z_1)=\left\{\frac{|\tau'|}{\pi}\ln 2+n\omega(\infty)\right\}K+iK'.
\end{equation}

Finally, with this notation we have the following expression for the complex
Green's function, see
\cite[\textsection 55, eq. (4)]{AKHef},
\begin{equation}\label{conf5}
\F(z)=\frac{\vartheta_1\left(\frac{u(z)+u(\infty)}{2K}\left|\right. -\frac 1\tau \right)}{\vartheta_1\left(\frac{u(z)-  u(\infty)}{2K}\left|\right. -\frac 1\tau \right)},
\end{equation}
where $\vartheta_1$ is the theta-function, which argument is shifted, comparably to $\vartheta_0$, by a half-period. To be precise,  see  Table VIII in \cite{AKHef},
$$
\vartheta_1(v+\frac \tau 2|\tau)=ie^{\frac{-\pi i\tau}{4}}e^{-\pi iv}\vartheta_0(v|\tau).
$$

Thus, due to  (\ref{ex2}) and (\ref{conf5}), we have
\bean
e^{g(z_1)}=|\Phi(z_1)|=\left|
\frac
{\vartheta_1\left(\frac{\{n\o(\infty)+|\tau'|\frac{\ln 2}{\pi} \}+\o(\infty)}{2}-\frac{1}{2\tau}\left|\right. -\frac 1\tau \right)}
{\vartheta_1\left(\frac{\{n\o(\infty)+|\tau'|\frac{\ln 2}{\pi} \}-\o(\infty)}{2}-\frac{1}{2\tau}\left|\right. -\frac 1\tau \right)}
\right|.
\eean
Using reduction by a half period we get
$$
e^{g(z_1)}=
\left|
\frac
{\vartheta_0\left(\frac{\{n\o(\infty)+|\tau'|\frac{\ln 2}{\pi} \}+\o(\infty)}{2}\left|\right. \tau' \right)}
{\vartheta_0\left(\frac{\{n\o(\infty)+|\tau'|\frac{\ln 2}{\pi} \}-\o(\infty)}{2}\left|\right. \tau' \right)}
\right|.
$$
By
$\mu(1^*,\G_n)=R(\i)e^{g(z_1)},$
in combination with the second expression in (\ref{ex7}), we obtain (\ref{ex}). \endproof

\bigskip

\centerline{******************************************}

\bigskip

The authors thank B. Nagy and S. Kalmykov
for bringing the paper \cite{TD} to their attention.

Vilmos Totik

MTA-SZTE Analysis and Stochastics Research Group

Bolyai Institute

University of Szeged

Szeged

Aradi v. tere 1, 6720, Hungary

\smallskip
and
\smallskip

Department of Mathematics and Statistics

University of South Florida

4202 E. Fowler Ave, CMC342

Tampa, FL 33620-5700, USA
\medskip

{\it totik@mail.usf.edu}
\vskip1cm

Peter Yuditskii

Johannes Kepler University Linz

Abteilung f¨ur Dynamische Systeme und Approximationstheorie

A-4040 Linz, Austria.
\medskip

{\it petro.yudytskiy@jku.at}

\begin{thebibliography}{9}
\bibitem{A1} N. I. Akhiezer,
\textit{\"Uber einige Funktionen welche in zwei gegebenen
Intervallen am wenigsten von Null abweichen, I--III},
{Izv. AN SSSR, 9 (1932), 1163--1202;}
4 (1933), 309--344, 499--536 (German).
Zbl 0007.00802, 0007.34106, 0007.34201. Russian translation:
N. I. Akhiezer, Selected works in Function theory and
mathematical physics, vol. 1, Akta, Kharkiv, 2001.
Zbl 1044.01014.

 \bibitem{AKHef} N.I. Akhiezer, \textit{Elements of the Theory of Elliptic Functions}. Amer. Math. Soc., Providence, RI, 1990.


\bibitem{Aptekarev1} A. I. Aptekarev, Asymptotic behavior of polynomials of simultaneous orthogonality and of a system of extremal problems for analytic functions, (Russian)
{\it Akad. Nauk SSSR Inst. Prikl. Mat.}, Preprint 1987, no. 168, 29 pp.

\bibitem{Aptekarev2} A. I. Aptekarev, Asymptotic properties of polynomials orthogonal on a system of contours, and periodic motions of Toda chains,  (Russian)
{\it Mat. Sb. (N.S.)}, {\bf 125}(167) (1984), no. 2, 231--258.

\bibitem{Faber} G. Faber, \"Uber Tschebyscheffsche Polynome, {\it J. Reine Angew. Math.},
{\bf  150}(1919), 79--106.


\bibitem{Gardiner}  D. H. Armitage and S. J. Gardiner,
{\it Classical Potential Theory}, Springer Verlag,
Berlin, Heidelberg, New York, 2002.

\bibitem{Garnett} J. B. Garnett and  D. E. Marshall,
{\it Harmonic measure}, Cambridge University Press,
New mathematical monographs,  Cambridge, New York, 2005.

\bibitem{Kal} V. A.  Kaliaguine, A note on the asymptotics of orthogonal polynomials on a complex arc: the case of a measure with a discrete part,
{\it  J. Approx. Theory}, {\bf  80}(1995),  138--145.






\bibitem{Peh1}
        {  F. Peherstorfer},
        Orthogonal and extremal polynomials on several
intervals,
        {\it J. Comp. Applied Math.}, {\bf 48}(1993), 187--205.


\bibitem{Peh2}
        {  F. Peherstorfer},
        Elliptic orthogonal and extremal polynomials,
        {\it J. London Math. Soc.}, {\bf 70}(1995), 605--624.

\bibitem{Peh3}
        {  F. Peherstorfer and K. Schiefermayr},
        Theoretical and numerical description of extremal
polynomials on
        several intervals I,
        {\it Acta Math. Hungar.}, {\bf 83}(1999), 27--58.

\bibitem{Peh4}
        {  F. Peherstorfer and K. Schiefermayr},
        Description of extremal polynomials on several intervals and their computation, I. II,
                {\it Acta Math. Hungar.}, {\bf 83}(1999), 71--102, 103--128.

\bibitem{Peh5}
        F. Peherstorfer,
        Deformation of minimizing polynomials and
approximation of several
        intervals by an inverse polynomial mapping,
        {\it J. Approx. Theory}, {\bf 111}(2001), 180--195.

\bibitem{Peh6} F. Peherstorfer and P. Yuditskii,
Asymptotic behavior of polynomials orthonormal on a homogeneous set,
{\it  J. Anal. Math.}, {\bf  89}(2003), 113--154.

\bibitem{Ransford} T. Ransford, {\it Potential Theory in the Complex
Plane}, Cambridge University Press,  Cambridge, 1995.


\bibitem{Sodin} M. L. Sodin and P. M. Yuditskii, Functions least deviating
from zero on closed subsets of the real axis, {\it Algebra i Analiz}, {\bf 4}(1992), 1-61;
English transl. in {\it St. Petersburg Math. J.}, {\bf 4}(1993), 201--249.

\bibitem{Schief} K. Schiefermayr, A lower bound for the minimum deviation of the
Chebyshev polynomial on a compact real set, {\it East J. Approx.},
{\bf 14}(2008), 65--75.



\bibitem{TD} J.-P. Thiran and C. Detaille, Chebyshev polynomials on circular
arcs and in the complex plane, {\it Prograss in Approximation Theory},
771--786, Academic Press, Boston, MA, 1991.

\bibitem{TotikCG} V. Totik, Chebyshev polynomials on compact sets,
{\it Potential Analysis}, 2013, DOI: 10.1007/s11118-013-9357-6.

\bibitem{Walsh} J. L. Walsh,
{\it Interpolation and Approximation by
Rational Functions in the Complex Domain}, third edition,
Amer. Math. Soc. Colloquium Publications, {\bf XX},
Amer. Math. Soc., Providence, 1960.


\bibitem{W} H. Widom, Extremal polynomials associated
with a system of curves in the complex plane, {\it
Adv. Math.}, {\bf 3}(1969), 127--232.
\end{thebibliography}
\end{document}